%% file: ex_article.tex
\documentclass[Draft,onefignum,onetabnum,a4paper]{siamart190516}
\usepackage[margin=1in,footskip=0.25in]{geometry}
\usepackage{float}
\input{ex_shared}

\ifpdf
\hypersetup{
  pdftitle={The mysteries of the best approximation and Chebyshev expansion for the function with logarithmic regularities},
  pdfauthor={Xiaolong Zhang}
}
\fi

\begin{document}

\maketitle

\begin{abstract}
  The best polynomial approximation and Chebyshev approximation are both important in numerical analysis. In tradition,
 the best approximation is regarded as more better than the Chebyshev approximation, because it is usually considered in
the uniform norm. However, it not always superior to the latter noticed by Trefethen \cite{Trefethen11sixmyths,Trefethen2020} for the algebraic singularity function. Recently Wang \cite{Wang2021best} have proved it in theory. In this paper, we find that for the functions with logarithmic regularities, the pointwise errors of Chebyshev approximation are smaller than the ones of the best approximations except only in the very narrow boundaries at the same degree. The pointwise error for Chebyshev series, truncated at the degree $n$ is $O(n^{-\kappa})$ ($\kappa = \min\{2\gamma+1, 2\delta + 1\}$), but is worse by one power of $n$ in narrow boundary layer near the weak singular endpoints. Theorems are given to explain this effect.
\end{abstract}



\section{Introduction}
\label{sec:intro}

The best polynomial approximations and the Chebyshev approximations are dense in numerical analysis and scientific computing with applications in rootfinding, function approximation, integral equations and differential equations. The best approximation is an old idea, which goes back to Chebyshev himself. Before the fast computers, the best approximations received more attentions than the alternatives. Thus in a long history, the best approximation is regarded as the optimal. At the early age of nineteen century, with the computer appearance and the Fast Fourier Transform being proposed \cite{CooleyTukey1965}, researchers have gradually begun to pay attentions to the Chebyshev approximation. In fact, the Chebyshev approximation is more practical and useful, because computing the best approximation is much more expensive. Furthermore, the Chebyshev approximation can easily be used to solve the quadrature and differential problems. During past several decades, it is widely used to solve many mathematical or engineering problems, for more details you can refer \cite{Boyd2001,Canuto2006,Guo1998,LiuWangLi2019,ShenTangWang2011,Trefethen2020,XiangLiu2020}. At the same time, connections and differences between best approximation and other orthogonal polynomial projection are discovered \cite{wang2020faster}. To begin, a bit more notation needs to be given.

In this paper, the interval $\Omega=[-1,1]$ is considered. $T_{k}(x) = \cos(k \arccos(x))$ is the Chebyshev polynomial of the first kind. If a function $f(x)$ satisfies the Dini-Lipschitz continuous condition, the following series can uniformly convergence to the function itself
\begin{equation}
f(x) = \sum_{k=0}^{\infty} {'}a_{k} T_{k}(x),\quad x\in [-1,1], \quad a_{k} = \frac{2}{\pi} \int_{-1}^{1} f(x)T_{k}(x)\frac{1}{\sqrt{1-x^2}} dx,
\end{equation}
where the prime on the sum denotes the first term is divided by 2. Truncate the exact series to a polynomial of degree $n$
\begin{equation}
f_{n} = \sum_{k=0}^{n} a_{k} T_{k}(x).
\end{equation}
\begin{figure}
\vspace{-0.5cm}
\includegraphics[scale=0.52]{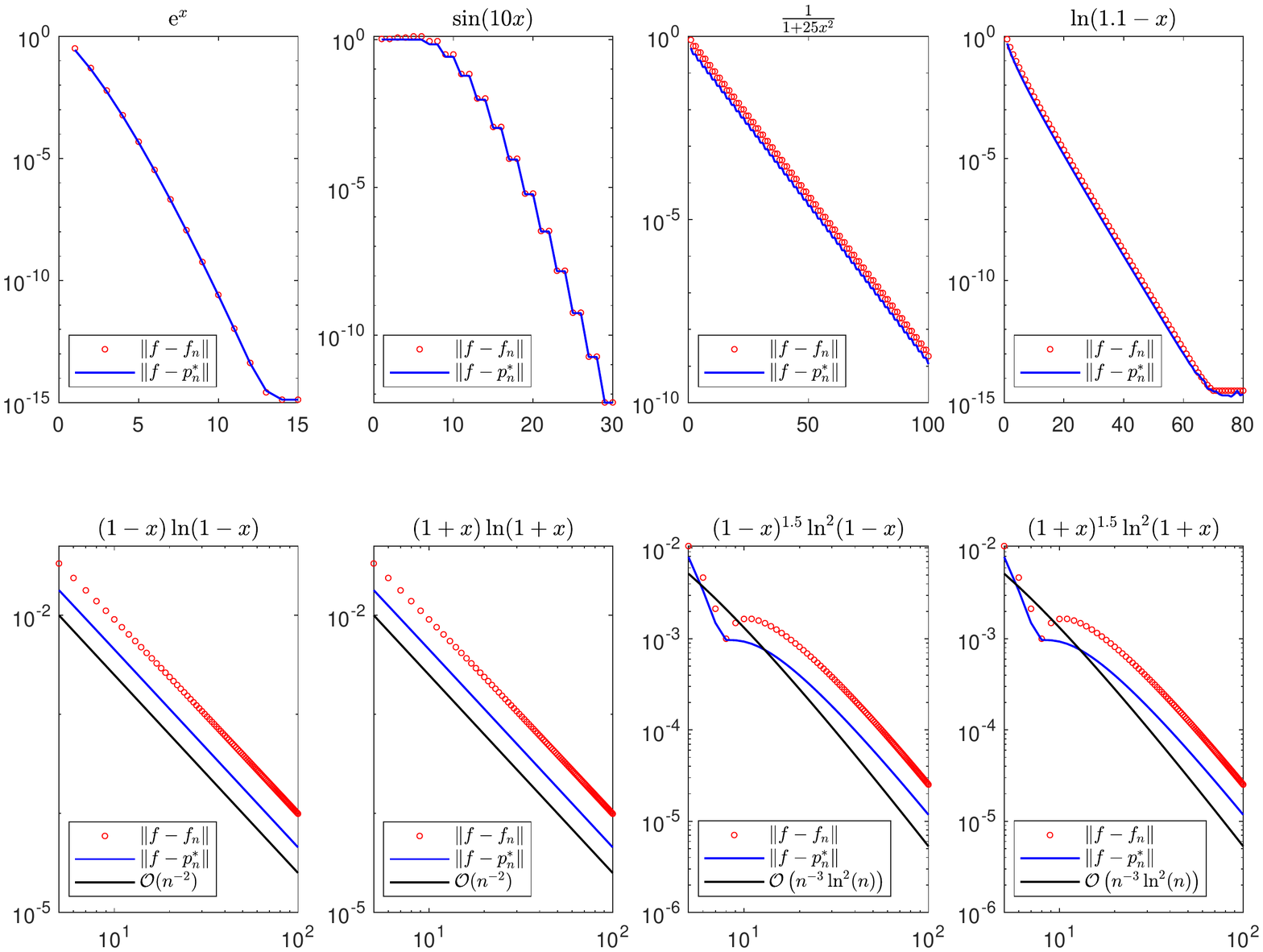}
\vspace{-1cm}
\caption{Comparisons between the errors of Chebyshev truncations and the errors of best approximation for the function $f(x)$, the horizon axis denotes the degree of approximations, the vertical axis denotes the error in the uniform norm. Red circle: error of Chebyshev truncation, blue solid line: error of best approximation. First row from the left to the right, the first : $f(x)=\mathrm{e}^{x}$, the second : $f(x) = \sin(10x)$, the third : $f(x) = \frac{1}{1+25x^2}$, the fourth : $f(x)=\ln(1.1-x)$; Second row from the left to the right, the first : $f(x) = (1-x)\ln(1-x)$, the second : $f(x) = (1+x)\ln(1+x)$, the third : $f(x) = (1-x)^{1.5}\ln^{2}(1-x)$, the fourth : $f(x) = (1+x)^{1.5}\ln^{2}(1+x)$.}
\label{fig:InfinityNormErrorsChebVsBest}
\end{figure}

Let $\mathcal{P}_n$ to be subspace of all real polynomial of degree at most $n$ in $C[-1,1]$. $p_{n}^{*}$ denotes the best polynomial approximation in $\mathcal{P}_{n}$ to function $f(x)$ in the maximum norm, i.e., $\parallel f - p_{n}^{*} \parallel_{L^{\infty}(\Omega)} = \min\limits_{h\in \mathcal{P}_{n}}\parallel f- h\parallel_{L^{\infty}(\Omega)}$.
Just as pointed out in \cite{Trefethen2020,Wang2021best}, the best approximation is usually more accurate than the same degree Chebyshev truncation in view of maximum norm as is shown in \cref{fig:InfinityNormErrorsChebVsBest}. The classical theorem can confirm this.
\begin{theorem}[Chebyshev truncation is near-best\cite{Trefethen2020}]\label{thm:ChebTrunNearBest}
Let f be continuous on $[-1,1]$ with degree $n$ Chebyshev truncation $f_n$ and best approximation $p_n^{*},\;n>1$. Then
\begin{equation}
\parallel f - f_n \parallel \le \left(4 + \frac{4}{\pi^{2}}\ln(n)\right) \parallel f - p^{*}_{n} \parallel.
\end{equation}
\end{theorem}

\begin{figure}[!ht]
\centering
\vspace{-0.5cm}
\includegraphics[scale=0.535]{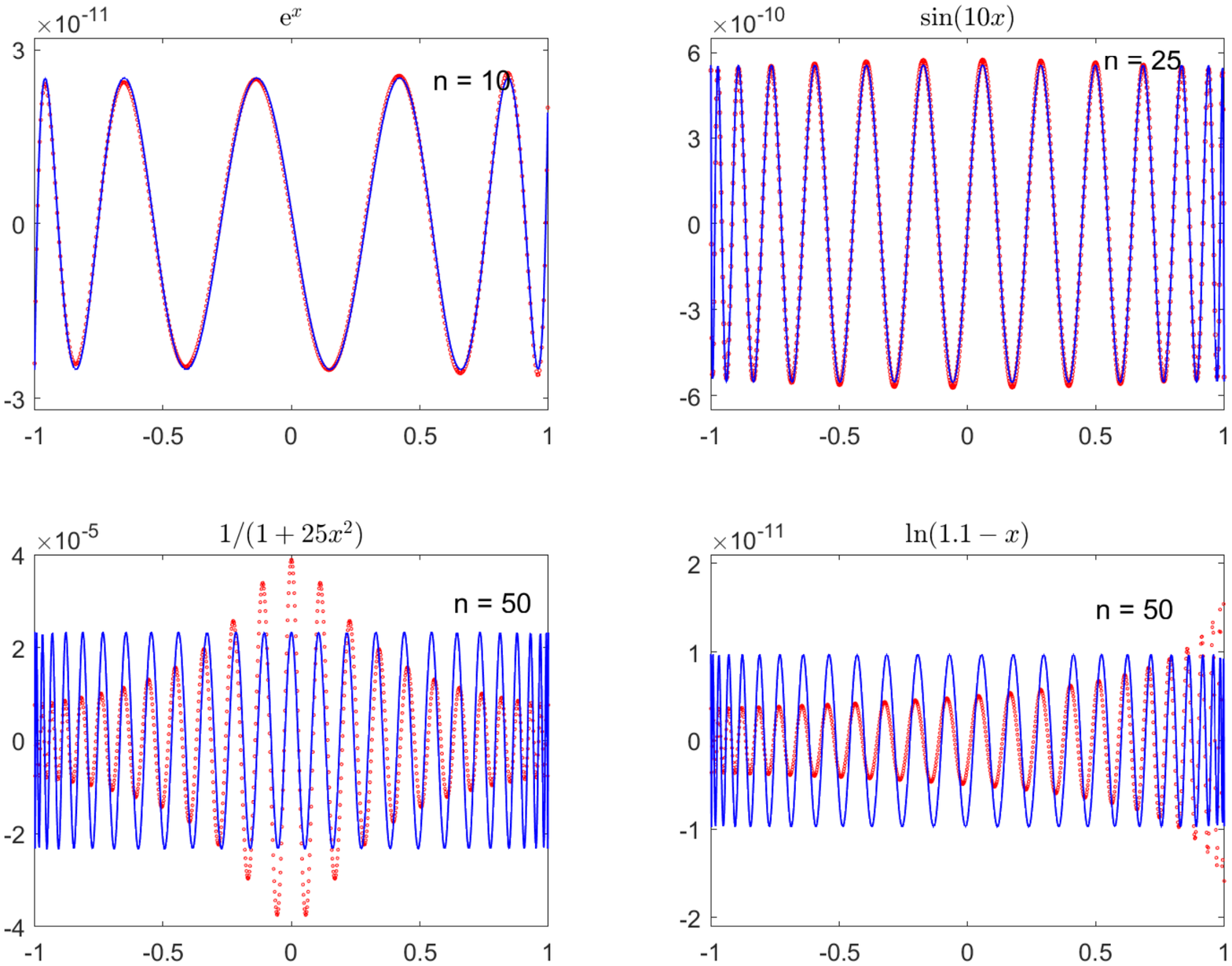}
\vspace{-1cm}
\caption{Comparisons of the Chebyshev approximations with the best approximations for the function $f(x),\;x\in[-1,1]$. Blue line: the pointwise
error of the best approximation of degree n, red circle: the pointwise error of the Chebyshev projection of degree n. Upper
left: $f(x) = \mathrm{e}^{x}$, Upper right: $f(x) = \sin(10\,x)$, Lower left: $f(x) = \frac{1}{1+25x^2}$, Lower right: $f(x) = \ln(1.1-x)$.}
\label{intro:ComparisonsGalleryFuns}
\end{figure}

\begin{figure}[!h]
\vspace{-3cm}
\centering
\includegraphics[scale=0.535]{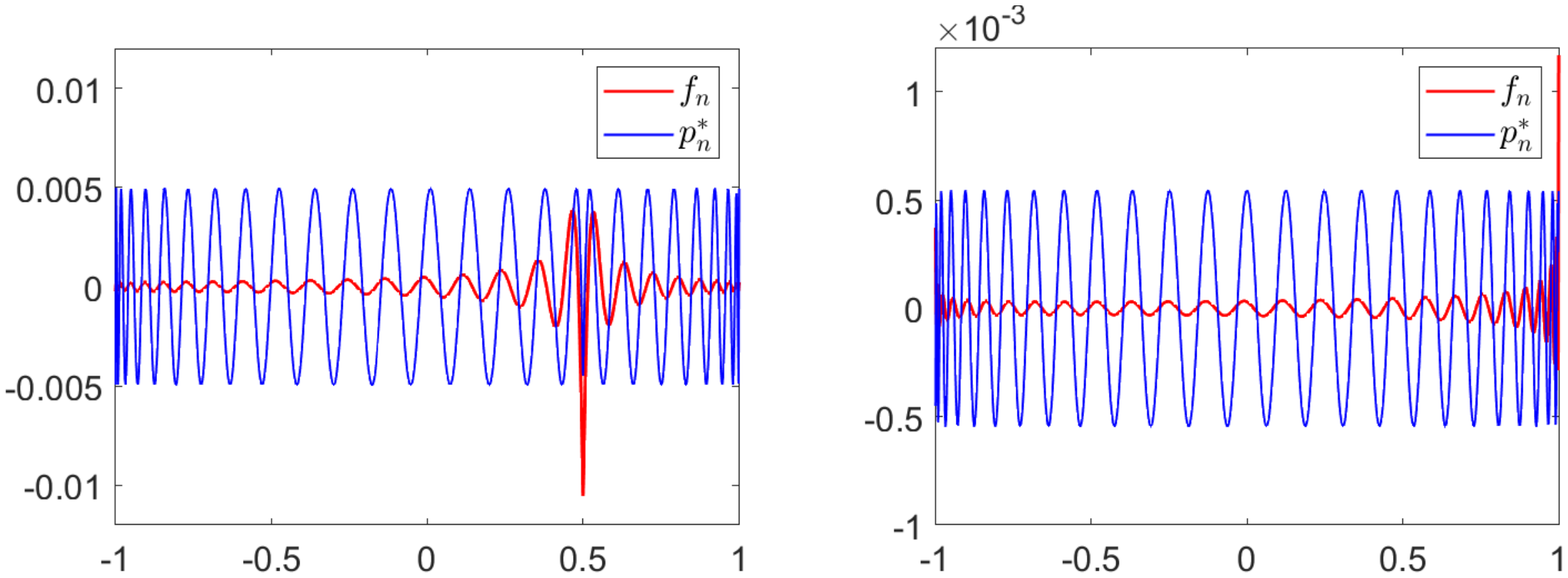}
\vspace{-3.5cm}
\caption{Pointwise Error curves of $f_n$ and $p_n^{*}$ for the function $f(x)$ with n = 50. Left: $f(x) = |x-0.5|$, Right:
$f(x) = \left(1+\frac{x}{2}\right)(1 - x^2)\ln(1 - x^2)$.}
\label{intro:AlgebraicVsLog}
\end{figure}

The proof can refer \cite{Trefethen2020}. Most people focus on the comparisons the two approximations in the uniform norm, but very little attention has been paid to the pointwise errors. Thus, on the face of it we seem to have the best approximation is always accurate more than its younger brother, the Chebyshev truncation. But recently, Trefethen \cite{Trefethen2020} have noticed the surprising phenomenon that the Chebyshev approximation is more accurate than the best approximation for the algebraic singular function
\begin{equation}
f(x) = |x - 0.25|^{\alpha}g(x), \quad x\in[-1,1],\quad \quad  \alpha >0,\quad g(x)\, \text{is analytic},
\end{equation}
without other singularities in the interval except in the narrow layer centered the algebraic singularity. Wang \cite{Wang2021best} gave an detail explanation in theory why the Chebyshev truncation $f_n$ is better than the best approximation $p_n^{*}$ except only in the very narrow layer centered at the singularity point. In fact, the behavior of the two approximations are rather bewildering as is shown in \cref{intro:ComparisonsGalleryFuns}. In this figure, all the functions has no singularities on the real interval $[-1,1]$. \cref{intro:AlgebraicVsLog} shows the pointwise error curves of the functions with singularities. By comparison the \cref{intro:ComparisonsGalleryFuns} and the \cref{intro:AlgebraicVsLog}, we can see that
\begin{itemize}
\item For the analytic function $f(x)$ on the complex plane, the both approximations seem almost equal that the difference hardly matters.
\item Although the function $f(x)$ has no singularities on the interval $[-1,1]$, but owns singularities on the whole complex plane. The pointwise error of the Chebyshev approximation near the singularities on the whole complex plane might be larger than the counterpart error of the best approximation.
\item Just like the algebraic singular function, the logarithm singular function owns a much better Chebyshev approximation rather than the best approximation for almost all values of $x\in [-1,1]$.
\end{itemize}

There are interesting phenomena as depicted in the right panel of \cref{intro:AlgebraicVsLog}, specifically,
\begin{enumerate}
\item [(i)] Why does the pointwise error of Chebyshev truncation blow up at the endpoints for the logarithm singularity function?
\item [(ii)] The both points $x = -1$ and $x = 1$ cause the weak singularity for the function
$f(x) = \left(1+\frac{x}{2}\right)(1 - x^2)\ln(1 - x^2)$. Why does the pointwise error of Chebyshev truncation at the point $x=-1$ is smaller than
 the counterpart error at the point $x = 1$?
\end{enumerate}

Inspired by the idea in the article \cite{Wang2021best}, we make above observations. To answer the questions, the pointwise error analysis of Chebyshev truncation is indispensable. In this paper, we will consider the following function
\begin{equation}\label{eq:orignalfun}
f(x) = (1 - x)^{\gamma}\ln^{\mu}(1 - x)g_{1}(x) + (1 + x)^{\delta}\ln^{\mu}(1 + x)g_{2}(x), \; x\in[-1,1],
\end{equation}
where $\gamma>-\frac{1}{2}, \delta>-\frac{1}{2}, g_{1}(x), g_{2}(x)$ are analytic on $ [-1,1], \mu$ is a positive integer. The function plays
an important role in function approximations and Frlory Huggins theory etc. Similar to the reference \cite{Wang2021best}, to explain the behavior of error for the logarithm singularity function, one needs to analyze the asymptotic behavior of the following functions
\begin{equation}\label{fun:PsiCos}
\Psi_{\nu}^{C}(t,n) = \sum_{k=n+1}^{\infty} \frac{\cos(k\,t)}{k^{\nu+1}},\quad \Psi_{\nu,\mu}^{C,L}(t,n) = \sum_{k=n+1}^{\infty}
\frac{\cos(k\,t)}{k^{\nu+1}}\ln^{\mu}(k)
\end{equation}
where $t\in[0,\pi], \nu>0, \mu \in \mathbb{N}$ . For the series $\Psi_{\nu}^{C}(t,n)$, Wang \cite{Wang2021best} has given the asymptotic representation as $n\to +\infty$. In \cref{sec:pointwiseChebTrun}, we will use a more simple technique to estimate the asymptotic behavior of $\Psi_{\nu}^{C}(t,n)$ and
 $\Psi_{\nu,\mu}^{C,L}(t,n)$.

The paper is organized as follows. In \cref{sec:aym}, the steepest descent method is used to provide the asymptotic coefficients of the Chebyshev
expansion for the function \cref{eq:orignalfun}. In \cref{sec:pointwiseChebTrun}, we discuss the pointwise errors of the Chebyshev truncations for the function \cref{eq:orignalfun}. In \cref{sec:interp}, based on the theory analyses in \cref{sec:pointwiseChebTrun}, the asymptotic behavior of pointwise error of Chebyshev interpolation is provided.

\section{Asymptotic Coefficients of Chebyshev expansions}
\label{sec:aym}

In this section, we will use the steepest descent method to give the asymptotic coefficients for the function \cref{eq:orignalfun}. The asymptotic estimate of Chebyshev coefficients for the function \cref{eq:orignalfun} when $\mu = 1$ was given by Boyd \cite{Boyd1989}. In this paper, the asymptotic coefficients for the function with the general positive integer $\mu =  1, 2,\cdots,$ are extended. Xiang \cite{Xiang2021} gave the Jacobi coefficients for the similar logarithmic singularity function. Now we consider the function
\begin{equation}\label{eq:funoneminius}
f(x) = (1- x )^{\gamma} \ln^{\mu}( 1 - x )g_{1}(x),\quad x\in[-1,1],
\end{equation}
where $\gamma>-\frac{1}{2}, \mu \in \mathbb{N}_{+}$ and $g_{1}(x)$ is analytic on $[-1,1]$. If $\mu =0, f(x)=(1-x)^{\gamma}$, the function may only
own branch singularity when $\gamma$ is not an integer. Before getting the asymptotic coefficients of Chebyshev expansion, two lemmas are given.

\begin{lemma}
Suppose $\mathrm{Re}(\gamma)>-1, \mathrm{Re}(k)>0$, then
\begin{equation}\label{eq:lemone}
\int_{0}^{\infty} s^{\gamma} \ln(s)\, \mathrm{e}^{-ks} ds  = \frac{1}{k^{\gamma+1}}\Gamma(\gamma+1)\left[\Psi(\gamma+1) - \ln(k)\right],
\end{equation}
where $\Psi(\gamma+1)$ is the digamma function.
\end{lemma}

\begin{proof}
Set $t = ks$, then the integral can be written as
\begin{equation*}
\begin{aligned}
\int_{0}^{\infty} s^{\gamma} \ln(s)\, \mathrm{e}^{-ks}ds& = \frac{1}{k^{\gamma+1}}\int_{0}^{\infty} t^{\gamma}[\ln(t) - \ln(k)] \mathrm{e}^{-t} dt \\
& = \frac{1}{k^{\gamma+1}} \Gamma(\gamma+1)\left(\Psi(\gamma+1) - \ln(k)\right).\\
\end{aligned}
\end{equation*}
\end{proof}

The lemma can also be found in identity 4.352 of Gradshteyn and Rhyzik \cite{Gradshteyn2014}, but the proof was not be provided, so the detail proof is given above.

\begin{lemma}\label{lem:two}
Suppose $\mathrm{Re}(\gamma)>-1, \mathrm{Re}(k)>0 $ and $\mu$ is a given positive integer, then
\begin{equation}
\begin{aligned}
\int_{0}^{\infty} s^{\gamma} \ln^{\mu}(s) \mathrm{e}^{-ks} ds &= \frac{1}{k^{\gamma+1}}
\sum_{\ell=0}^{\mu}\binom{\mu}{\ell}\Gamma^{(\mu-\ell)}(\gamma+1)(-1)^{\ell}\ln^{\ell}(k)\\
&\sim (-1)^{\mu}\Gamma(\gamma+1) \frac{\ln^{\mu}(k)}{k^{\gamma+1}}, \quad \text{as}\; n\to +\infty.
\end{aligned}
\end{equation}
\end{lemma}
\begin{proof}
To proof this lemma, one only needs to take the $(\mu -1)$th derivative of Eq.\cref{eq:lemone} about $\gamma$.

\end{proof}

\begin{theorem}\label{thm:asycoeffsfunminus}
If a function $f(x)$ has a singularity at $x = 1$ of the form \cref{eq:funoneminius}, then the coefficients of Chebyshev series are asymptotically
given by
\begin{equation}
\begin{aligned}
a_{k} \sim &\frac{2^{\mu + 1 - \gamma}}{\pi}\, g_{1}(1)\,(-1)^{\mu+1}\,\sin(\gamma\,\pi)\,\Gamma(2\gamma+1)\,\frac{\ln^{\mu}(k)}{k^{2\gamma+1}}+
 (-1)^{\mu}\,\mu\, 2^{\mu-\gamma}\,  g_{1}(1)\\
& \times\left( \frac{2}{\pi}\sin(\gamma\,\pi)\Gamma'(2\gamma+1) +  \left(\cos(\gamma\,\pi)-\frac{\ln(2)}{\pi}\right)\Gamma(2\gamma+1)\right)
\frac{\ln^{\mu-1}(k)}{k^{2\gamma+1}} + \mathcal{O}\left(\frac{\ln^{\mu-2}(k)}{k^{2\gamma + 1}}\right).
\end{aligned}
\end{equation}
\end{theorem}

\begin{proof}
The Chebyshev expansion can be expressed as
\begin{equation*}
f(x) = \sum_{k=0}^{\infty} {'} a_k T_k(x),
\end{equation*}
where
\begin{equation*}
a_{k} = \frac{2}{\pi} \int_{-1}^{1} f(x) \frac{1}{\sqrt{1-x^2}}T_{k}(x)dx,\quad k>0.
\end{equation*}
Set $x = \cos(t), t\in[0,\pi]$, the coefficients $a_{k}$ are
\begin{equation*}
a_{k} = \frac{2}{\pi}\int_{0}^{\pi} f(\cos(t))\cos(k\,t) dt, \quad k\in \mathbb{N}_{+}.
\end{equation*}
Since the Chebyshev Coefficients are dominated by the worst singularities, thus we only need to evaluate the worst singularity integral.
Consequently, the asymptotic coefficients can be simplified as
\begin{equation*}
\begin{aligned}
a_{k} &\sim \frac{2^{\mu+1-\gamma}}{\pi}g_{1}(1) \int_{0}^{\pi} t^{2\gamma} \ln^{\mu}(t) \cos(k\,t)dt\\
&\sim \frac{2^{\mu+1-\gamma}}{\pi}g_{1}(1) \mathrm{Re} \left\{\int_{0}^{\pi} t^{2\gamma} \ln^{\mu}(t)\, \mathrm{e}^{ikt} dt\right\},\quad k\ge1.
\end{aligned}
\end{equation*}
The special case $\gamma =0 , \mu =1$ are used to illustrate the steepest descent method \cite{Bender1978}. The key idea is to replace the $\cos(nt)$ by
$\mathrm{Re}(\mathrm{e}^{int})$ and deform the contour of integration into three line segments:
\begin{enumerate}
\item [(a)] $t=0 + i\cdot0$ to $t =0+  i\cdot\infty$.
\item [(b)] $t =0+ i\cdot \infty $ to $t = \pi+ i\cdot\infty$.
\item [(c)] $t = \pi + i\cdot\infty $ to $t = \pi + i\cdot 0$
\end{enumerate}
The second contour integration approximates zero. The third is much less than $\mathcal{O}(k^{-2\gamma-1})$. Thus the main contribution comes
from the integration in the first segment. The asymptotic expression of the first is
\begin{equation*}
a_{k} \sim \frac{2^{\mu+1-\gamma}}{\pi}g_{1}(1) \mathrm{Re} \left\{\int_{0+i0}^{0+i\infty} t^{2\gamma} \ln^{\mu}(t)\, \mathrm{e}^{ikt} dt\right\}
\end{equation*}
Set $t = is$ and
\begin{equation*}
I_{1} = \int_{0+i0}^{0+i\infty} t^{2\gamma} \ln^{\mu}(t)\, \mathrm{e}^{ikt} dt ,
\end{equation*}
thus, the integration becomes
\begin{equation*}
I_{1} = (i)^{2\gamma+1}\int_{0}^{\infty} s^{2\gamma} \ln^{\mu}(is)\, \mathrm{e}^{-ks} ds
\end{equation*}
Applying the identity $\ln(is) = \ln(s) + \frac{\pi}{2}i$ and taking the dominant terms, then the integration $I_{1}$  simplifies as
\begin{equation}
\begin{aligned}
I_{1} &\sim  (i)^{2\gamma+1} \int_{0}^{\infty} s^{2\gamma} \sum_{\ell =0}^{\mu}\binom{\mu}{\ell}
\ln^{\mu -\ell}(s)\left( \frac{\pi}{2}i\right)^{\ell}  \\
& \sim i^{2\gamma+1} \int_{0}^{\infty} s^{2\gamma} \ln^{\mu}(s)\mathrm{e}^{-ks}ds +\frac{\pi}{2} \mu \,i^{2\gamma+2}
 \int_{0}^{\infty} s^{2\gamma} \ln^{\mu-1}(s)\mathrm{e}^{-ks}ds+\cdots.
\end{aligned}
\end{equation}
By \cref{lem:two}, the asymptotic coefficients of Chebyshev expansion are finally obtained.
\end{proof}

\begin{figure}[!ht]
\centering
\vspace{-0.5cm}
\includegraphics[scale=0.535]{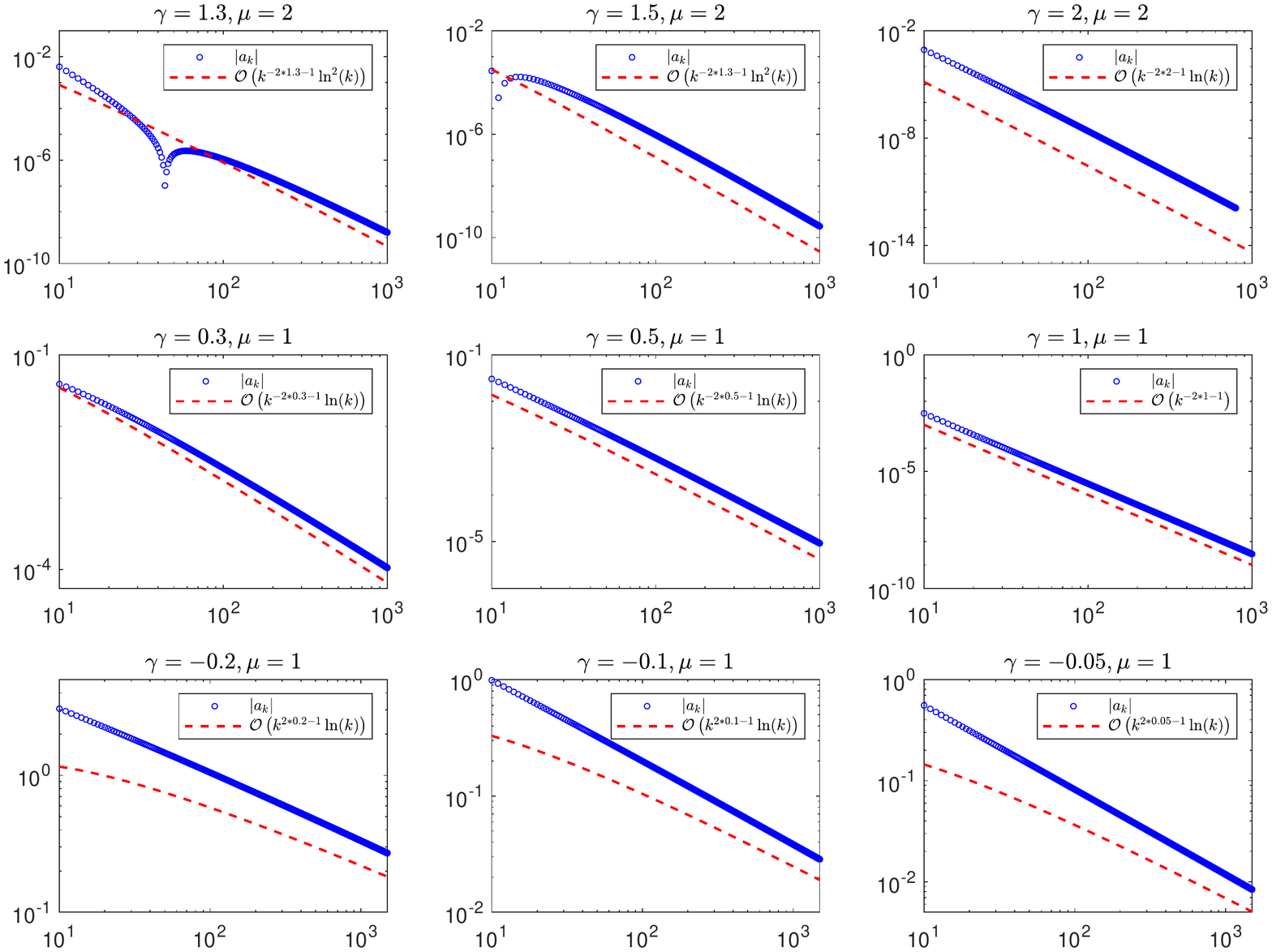}
\vspace{-1cm}
\caption{The Chebyshev coefficients for the function $f(x) = \left(1+\frac{x}{2}\right)(1-x)^{\gamma}\ln^{\mu}(1-x)$ with different $\gamma$
and $\mu$.}
\label{fig:funminuscoeffs}
\end{figure}

It is easy to see that when $\gamma$ is an integer in \cref{eq:funoneminius}, the order of $\ln(n)$ in coefficients will reduce one as is shown
in \cref{fig:funminuscoeffs}. It is also shown in \cite{Xiang2021}. In \cite{Boyd1989}, Boyd pointed out the special case when $\gamma$ is an integer and $\mu =1$, the $\ln(n)$ will disappear, which is consistent with the current result.

Next consider the Chebyshev coefficients for the function
\begin{equation}\label{eq:funtwoplus}
f(x) = (1 + x )^{\delta} \ln^{\mu} (1 + x)g_2(x), \quad x\in[-1,1],
\end{equation}
where $\delta>0, \mu\in \mathbb{N}_{+}$ and $g_{2}(x)$ is analytic on the interval $[-1,1]$. Here the asymptotic Chebyshev coefficients are given
in the following.
\begin{corollary}\label{Coro:CoeffsPLUS}
If a function $f(x)$ has a singularity at $x = -1$ of the form \cref{eq:funtwoplus}, then the coefficients
of Chebyshev series are asymptotically proportion to
\begin{equation}
\begin{aligned}
a_{k} \sim &\frac{2^{\mu + 1 - \delta}}{\pi}\, g_{2}(-1)\,(-1)^{k+\mu+1}\,\sin(\delta\,\pi)\,\Gamma(2\delta+1)\,
\frac{\ln^{\mu}(k)}{k^{2\delta+1}}+ (-1)^{k+\mu}\,\mu\, 2^{\mu-\delta}\,  g_{2}(-1)\\
& \times\left( \frac{2}{\pi}\sin(\delta\,\pi)\Gamma'(2\delta+1) +  \left(\cos(\delta\,\pi)-\frac{\ln(2)}{\pi}\right)\Gamma(2\delta+1)\right)
\frac{\ln^{\mu-1}(k)}{k^{2\delta+1}} + \mathcal{O}\left(\frac{\ln^{\mu-2}(k)}{k^{2\delta + 1}}\right).
\end{aligned}
\end{equation}
\end{corollary}
\begin{proof}
By the \cref{thm:asycoeffsfunminus}, the asymptotic Chebyshev coefficients $\tilde{a}_{k}$ for the function $f(-x)$ is
\begin{equation}
\begin{aligned}
\tilde{a}_{k} \sim &\frac{2^{\mu + 1 - \delta}}{\pi}\, g_{2}(-1)\,(-1)^{\mu+1}\,\sin(\delta\,\pi)\,\Gamma(2\delta+1)\,
\frac{\ln^{\mu}(k)}{k^{2\delta+1}}+ (-1)^{\mu}\,\mu\, 2^{\mu-\delta}\,  g_{2}(-1)\\
& \times\left( \frac{2}{\pi}\sin(\delta\,\pi)\Gamma'(2\delta+1) +  \left(\cos(\delta\,\pi)-\frac{\ln(2)}{\pi}\right)\Gamma(2\delta+1)\right)
 \frac{\ln^{\mu-1}(k)}{k^{2\delta+1}} + \mathcal{O}\left(\frac{\ln^{\mu-2}(k)}{k^{2\delta + 1}}\right).
\end{aligned}
\end{equation}
Since the symmetries of the Chebyshev polynomials $T_k(-x) = (-1)^{k} T_{k}(x)$, it is not difficult to find $a_{k} = (-1)^{k} \tilde{a}_{k}$.
\end{proof}
\begin{figure}[!h]
\vspace{-0.3cm}
\centering
\includegraphics[scale=0.535]{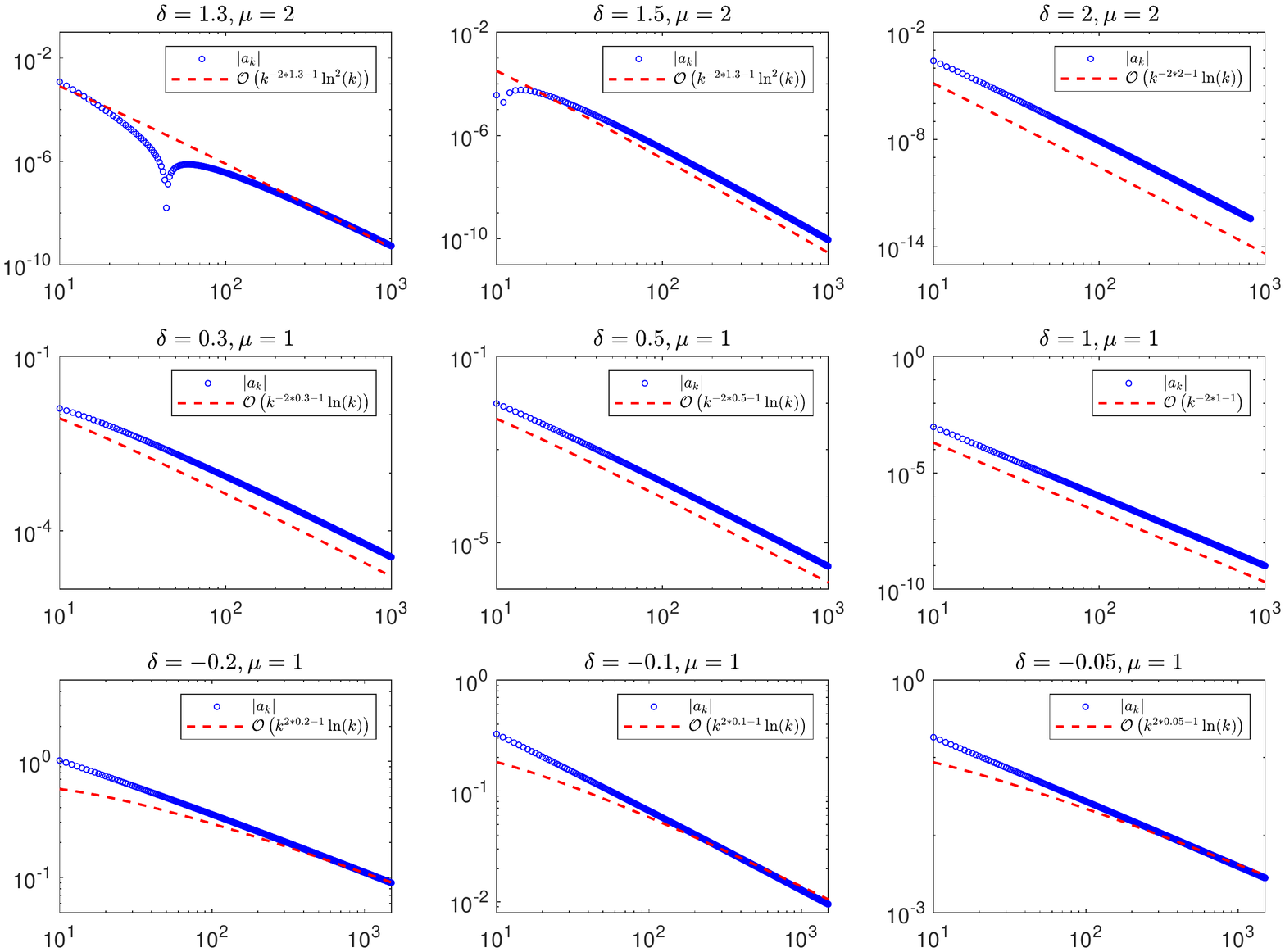}
\vspace{-1cm}
\caption{The Chebyshev coefficients for the function $f(x) = \left(1+\frac{x}{2}\right)(1+x)^{\gamma}\ln^{\mu}(1+x)$ with different $\gamma$
and $\mu$.}
\label{fig:CoeffsPlus}
\end{figure}

Before this, Boyd \cite{Boyd1989} have noticed the similar findings for $\mu = 1$. To illustrate the rules of asymptotic coefficient given in
 \cref{Coro:CoeffsPLUS}, the coefficients for the function $f(x) = (1+\frac{x}{2})(1+x)\ln^{\mu}(1+x)$ are shown in \cref{fig:CoeffsPlus}.
Then we consider more general weak singularity function \cref{eq:orignalfun}.

\begin{corollary}\label{coro:coeffs}
If the function \cref{eq:orignalfun} has weak singularities at $ x = 1, x=-1$, the asymptotic Chebyshev coefficients are
\begin{equation}
|a_{k}| =
\begin{cases}
\mathcal{O}\left( \frac{ln^{\mu}(k)}{k^{\min\{2\gamma+1, 2\delta+1\}}}\right), \quad \gamma\notin \mathbb{Z}_{+},\; \delta\notin \mathbb{Z}_{+},\\
\mathcal{O}\left( \frac{ln^{\mu-1}(k)}{k^{\min\{2\gamma+1, 2\delta+1\}}}\right), \quad \gamma\in\mathbb{Z}_{+},\;\delta\in\mathbb{Z}_{+}\\
 \max\left\{\mathcal{O}\left(\frac{\ln^{\mu-1}(k)}{k^{2\gamma+1}}\right),\mathcal{O}\left(\frac{\ln^{\mu}(k)}{k^{2\delta+1}}\right)\right\},\quad
\gamma\in \mathbb{Z}_{+},\; \delta\notin \mathbb{Z}_{+},\\
 \max\left\{\mathcal{O}\left(\frac{\ln^{\mu}(k)}{k^{2\gamma+1}}\right),\mathcal{O}\left(\frac{\ln^{\mu-1}(k)}{k^{2\delta+1}}\right)\right\},
\quad \gamma\notin \mathbb{Z}_{+},\; \delta\in \mathbb{Z}_{+},.\\
\end{cases}
\end{equation}
\end{corollary}

\begin{figure}[!h]
\centering
\includegraphics[height=3.5in,width=6.2in]{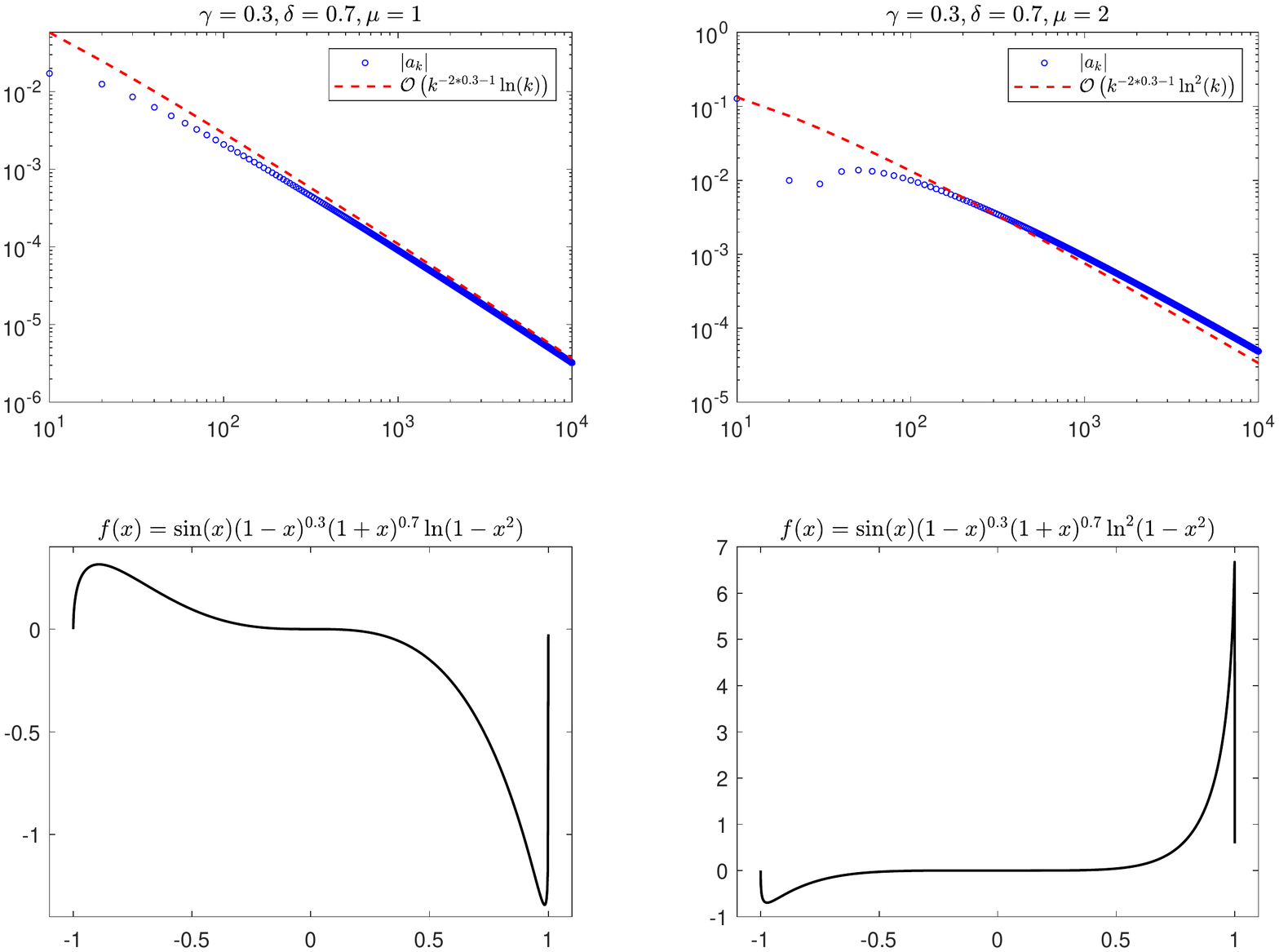}
\vspace{-0.5cm}
\caption{The Cheybshev coefficients for the function $f(x) = \sin(x)(1-x)^{\gamma}(1+x)^{\delta}\ln(1-x^2) $ \cite{Xiang2021}. Upper: Chebyshve coefficients;
 Lower: the curves of the functions. }\label{fig:fun1plusminusCoeffs}
\end{figure}

The Chebyshev coefficients decay rate for the function $f(x) = \sin(x)(1-x)^{\gamma}(1+x)^{\delta}\ln(1-x^2)$ are illustrated in
\cref{fig:fun1plusminusCoeffs}. It is shown that a very narrow peak near the right boundary layer appears in the curves of the functions.
When the $\mu$ increases larger, the peak becomes more narrow and shaper. That is the reason why the Chebyshev coefficients for the function
 $f(x)$ with $\mu=2$ decay slower.

\begin{remark}
In \cref{coro:coeffs}, the asymptotic behavior also proved by Xiang \cite{Xiang2021}, different with his method, the steepest descent method is used in this article.
\end{remark}

\section{Pointwise error estimates of Chebyshev truncation}
\label{sec:pointwiseChebTrun}

In this section, we will investigate the pointwise error estimates of Chebyshev truncation for the function \cref{eq:orignalfun}.
Before this, we give two lemmas.

\begin{lemma}\label{Lem:PsiAsym}
Let $\Psi_{\nu}^{C}(t,n)$ be the function defined in \cref{fun:PsiCos}. For $n\gg 1$, the following statement holds.
\begin{enumerate}
\item[(i)] If $t(\mathrm{mod}, 2\pi)\neq 0$ and $\nu>0$, then
\begin{equation}\label{eq:PsiAsym}
\Psi_{\nu}^{C}(t,n) = -D_{n}(t)\,n^{-\nu-1} + \mathcal{O}(n^{-\nu-2}),
\end{equation}
where $D_n(t)=\frac{\sin\left(\frac{2n+1}{2}t\right)}{2\sin(\frac{t}{2})}$ is Dirichlet's kernel.
\item[(ii)] If $t(\mathrm{mod}, 2\pi)= 0$ and $\nu>0$, then
\begin{equation}\label{eq:PsiAsym2}
\Psi_{\nu}^{C}(t,n) = \frac{1}{\nu} n^{-\nu}-\frac{1}{2}n^{-\nu-1} + \mathcal{O}(n^{-\nu-2}).
\end{equation}
\end{enumerate}
\end{lemma}
\begin{proof}
The \cref{Lem:PsiAsym} is proved by Wang \cite{Wang2021best} used the asymptotic presentation of the Hurwitz zeta function.
\end{proof}

Here an alternative method is provided to give a sharp upper bound of the $\Psi_{\nu}^{C}(x,n)$. For $x\in[\delta_{1},2\pi-\delta_{1}]$, $0<\delta_{1}\ll 2\pi$ is a small positive parameter, one has
\begin{equation*}
\left|\sum_{k = n+1}^{n+p} \cos(kt)\right| = \left|\frac{\sin\left(\frac{2n+2p+1}{2}t\right)- \sin\left(\frac{2n+1}{2}t\right) }
{2\sin(\frac{t}{2})}\right|\le \frac{1}{|\sin(\frac{t}{2})|},\quad t\in[\delta_{1},2\pi-\delta_{1}]
\end{equation*}
and the series $\{\frac{1}{k^{\nu}}\}$ is monotonic and converges to zero as $k\to \infty$. By the Dirichlet test theorem, one can conclude that
the series
\begin{equation*}
\Psi_{\nu}^{C}(t,n) = \sum_{k=n+1}^{\infty} \frac{\cos(k\,t)}{k^{\nu+1}}
\end{equation*}
is uniformly convergent on $t\in[\delta_{1}, 2\pi-\delta_{1}]$. Thus, we have
\begin{equation}\label{eq:PsiTrun}
\left|\sum_{k=n+1}^{n+p} \frac{\cos(k\,t)}{k^{\nu+1}}\right|\le \frac{2}{\sin(\frac{t}{2})}\left|\frac{1}{(n+1)^{\nu+1}}+2
\frac{1}{(n+p)^{\nu+1}}\right|,\quad p\in \mathbb{N}_{+}.
\end{equation}
By taking the limit $p\to \infty$ in \cref{eq:PsiTrun}, it holds, finally
\begin{equation*}
\left|\Psi_{\nu}^{C}(t,n)\right| \le \frac{2}{\sin(\frac{t}{2})}\frac{1}{(n+1)^{\nu+1}}= \mathcal{O}(n^{-\nu-1}), \quad n\gg 1.
\end{equation*}

\begin{remark}
In \cite{Wang2021best}, it is required $\nu>0$, but in current estimate we extend the interval to $\nu>-1$.
\end{remark}

\begin{lemma}\label{Lem:PsiLogAsym}
Let the function $\Psi_{\nu,\mu}^{C,L}(t,n)$ is defined in \cref{fun:PsiCos}, and
\begin{enumerate}
\item [(i)]
if $t (\mathrm{mod}, 2\pi)\neq 0$ and $\nu>0$, then there exists the asymptotic expression
\begin{equation}\label{eq:AsymPsiCosLn}
\Psi_{\nu,\mu}^{C,L}(t,n) = -D_{n}(t)\,n^{-\nu-1}\ln^{\mu}(n) + \mathcal{O}\left(n^{-\nu-2}\ln^{\mu}(n)\right).
\end{equation}
\item [(ii)]
if $t (\mathrm{mod}, 2\pi)= 0$ and $\nu>0$, then there exists the asymptotic expression
\begin{equation}\label{eq:AsymPsiCosLnBoundary}
\Psi_{\nu,\mu}^{C,L}(t,n) \sim \frac{1}{\nu} n^{-\nu}\ln^{\mu}(n)-\frac{1}{2} n^{-\nu-1} \ln^{\mu}(n) + \mathcal{O}\left(n^{-\nu-2}\ln^{\mu}(n)\right).
\end{equation}
\end{enumerate}
\end{lemma}

\begin{proof}
We will prove (i),(ii) respectively in the following.
\begin{enumerate}
\item[(i)]
For any $t\in \mathbb{R}$, if $t(\mathrm{mod},2\pi)\neq 0$, there must exist a positive small parameter satisfies $\delta_{1} \ll 2\pi$ such that
the series
\begin{equation*}
 \Psi_{\nu,\mu}^{C,L}=\sum\limits_{k=n+1}^{\infty} \frac{\cos(k\,t)}{k^{\nu+1}}\ln^{\mu}(k)
\end{equation*}
is uniformly convergent on the interval $[\delta_{1}, 2\pi-\delta_{1}]$. Thus, taking the derivative about $\nu$ in both sides of \cref{eq:PsiAsym} $\mu$ times, gives
\begin{equation*}
\Psi_{\nu,\mu}^{C,L}(t,n) = -\frac{\sin\left(\frac{2n+1}{2}t\right)}{2\sin(\frac{t}{2})}n^{-\nu-1}\ln^{\mu}(n) + \mathcal{O}\left(n^{-\nu-2}\ln^{\mu}(n)\right).
\end{equation*}
\item [(ii)] To prove this, we only need compute the $\mu$th derivative about $\nu$ in both sides of \cref{eq:PsiAsym2}.
\end{enumerate}

\end{proof}

Similar to the analysis above, for (i), we can also use the Dirichlet test theorem to estimate an upper bound of the function
$\Psi_{\nu,\mu}^{C,L}(x,n)$, satisfying
\begin{equation*}
\left|\Psi_{\nu,\mu}^{C,L}(t,n)\right|\le \frac{2}{\sin(\frac{t}{2})}\left(n^{-\nu-1}\ln^{\mu}(n)\right)= \mathcal{O}(n^{-\nu-1}\ln^{\mu}(n)).
\end{equation*}


A detailed understanding of \cref{Lem:PsiLogAsym} will lead insight into the pointwise error estimates.
Here we state the main results about pointwise error estimates as \cref{thm:errests}.

\begin{theorem}\label{thm:errests}
Let $f$ be the function defined in \cref{eq:funoneminius} and $f_{n}$ be Chebyshev truncation of degree $n$. For $n\to \infty$, the following
establishes
\begin{enumerate}
\item[(i)] if $x\in[-1,1)$, then
\begin{equation*}
f-f_{n} \sim \frac{\alpha_{1}}{\pi}\, \sin(\gamma\pi)\, \Gamma(2\gamma + 1)\,U_{2n}\left(\frac{x}{2}\right)\,
\frac{ \ln^{\mu}(n) }{ n^{ 2\gamma + 1 } } \\ - \frac{1}{2}\,\mu\, \alpha_{1}\,\alpha_{2}\, U_{2n}\left(\frac{x}{2}\right)
\frac{ \ln^{\mu-1 }(n) }{ n^{ 2\gamma +1} } + \mathcal{O}\left(\frac{ \ln^{\mu-2}(n) } { n^{2\gamma + 1} }\right).
\end{equation*}
where
\begin{displaymath}
\begin{aligned}
&\alpha_{1} = (-1)^{\mu}\,2^{\mu-\gamma}\,g_{1}(1)\\
&\alpha_{2} = \frac{2}{\pi}\,\sin(\gamma\pi)\, \Gamma'(2\gamma + 1) + \left( \cos(\gamma\pi)-\frac{\ln(2)}{\pi} \right)\Gamma( 2\gamma + 1).
\end{aligned}
\end{displaymath}
\item [(ii)]if $x=1$, then
\begin{equation*}
f-f_{n} \sim \frac{1}{\gamma}\, \frac{\alpha_{1}}{\pi}\, \Gamma( 2\gamma + 1 )\, \sin(\gamma \pi)\, \frac{ \ln^{\mu}(n) }{ n^{2\gamma } }  +
\frac{\mu}{2 \gamma} \, \alpha_{1} \, \alpha_{2} \, \frac{ \ln^{\mu-1}(n) }{ n^{2\gamma }} +
\mathcal{O}\left( \frac{ \ln^{\mu-2}(n) }{ n^{2\gamma} } \right)
\end{equation*}
with $\alpha_{1}, \alpha_{2}$ is given in (i) above.
\end{enumerate}
\end{theorem}
\begin{proof}
For the function \cref{eq:funoneminius}, the error of the Chebyshev truncation
\begin{equation*}
f(x)- f_{n} = \sum_{k = n+1}^{\infty} a_{k} T_{k}(x)= a_{k} \cos(k\,t),\quad n\gg 1,\quad x\in[-1,1],\quad x = \cos(t),
\end{equation*}
is considered. First, we prove the $(i)$, by the \cref{thm:asycoeffsfunminus}, the error can be written as
\begin{equation*}
f(x) - f_{n} \sim \sum_{k = n+1}^{\infty}\left(-\frac{2}{\pi} \, \alpha_{1} \, \sin(\gamma \, \pi) \, \Gamma( 2\gamma + 1 ) \,
\frac{ \ln^{\mu}(k) }{ k^{ 2\gamma + 1 }} + \mu\,\alpha_{1}\,\alpha_{2}\, \frac{\ln^{ \mu-1 }(k)}{ k^{ 2\gamma +1 }}\right)\, T_{k}(x) +
\mathcal{O}\left( \frac{ \ln^{ \mu - 2 }(n)}{ n^{ 2\gamma + 1 }}\right).
\end{equation*}
Due to $x\in[-1,1)$, which corresponds the variable $t\in(0,\pi]$. Combining this and the identity \cref{eq:AsymPsiCosLn} given in
\cref{Lem:PsiAsym}, establishes the item $(i)$. Next, taking into account that $x=1$, by the identity $T_{k}(1) =1$, the error gives
\begin{equation*}
f(x) - f_{n} \sim \sum_{k = n+1}^{\infty}\left(-\frac{2}{\pi} \, \alpha_{1} \, \sin(\gamma \, \pi)\, \Gamma( 2\gamma + 1 ) \,
\frac{ \ln^{\mu}(k) }{ k^{ 2\gamma + 1 }} + \mu\,\alpha_{1}\,\alpha_{2}\,\frac{\ln^{ \mu-1 }(k)}{ k^{ 2\gamma +1 }}\right) +
\mathcal{O}\left(\frac{ \ln^{\mu-2} } { n^{2\gamma} }\right),\quad n\gg 1.
\end{equation*}
Using the asymptotic expression \cref{eq:AsymPsiCosLnBoundary}, the item (ii) is proved.
\end{proof}

It turns out that the pointwise error at $x=1$ is $O(n)$ larger than the error at other point in $x\in[-1,1)$, a fact that plays an important role in approximation theory. The main reason is it exists a weak singularity at $x=1$.

\begin{theorem}
Let $f$ be the function defined in \cref{eq:funtwoplus} and $f_{n}$ be Chebyshev truncation of degree $n$. For $n\to \infty$, the following
statements establishes
\begin{enumerate}
\item[(i)] if $x\in(-1,1]$, then

\begin{equation*}
\hspace{-0.15in}
f - f_{n} \sim \frac{\beta_{1}}{\pi} \sin(\delta\,\pi) \Gamma(2\delta + 1)U_{2n}\left(\frac{x + \pi}{2}\right)
\frac{ \ln^{\mu}(n) }{ n^{ 2\delta + 1 } } - \frac{\mu}{2} \beta_{1}\beta_{2} U_{2n}\left(\frac{x+\pi}{2}\right)
\frac{ \ln^{\mu-1 }(n) }{ n^{ 2\delta +1} } + \mathcal{O}\left(\frac{ \ln^{\mu-2}(n) } { n^{2\delta + 1} }\right).
\end{equation*}
where
\begin{displaymath}
\begin{aligned}
&\beta_{1} = (-1)^{\mu}\,2^{\mu-\delta}\,g_{2}(-1)\\
&\beta_{2} = \frac{2}{\pi}\,\sin(\delta\pi)\, \Gamma'(2\delta + 1) + \left( \cos(\delta\pi)-\frac{\ln(2)}{\pi} \right)\Gamma( 2\delta + 1).
\end{aligned}
\end{displaymath}
\item [(ii)]if $x=-1$, then
\begin{equation*}
f-f_{n} \sim \frac{1}{\delta}\, \frac{\beta_{1}}{\pi}\, \Gamma( 2\delta + 1 )\, \sin(\delta\, \pi)\, \frac{ \ln^{\mu}(n) }{ n^{2\delta } }  +
\frac{\mu}{2 \delta} \, \beta_{1} \, \beta_{2} \, \frac{ \ln^{\mu-1}(n) }{ n^{2\gamma }} +
\mathcal{O}\left( \frac{ \ln^{\mu-2}(n) }{ n^{2\delta} } \right)
\end{equation*}
with $\beta_{1}, \beta_{2}$ is given in (i) above.
\end{enumerate}
\end{theorem}
\begin{proof}
To prove this theorem, the identity
\begin{equation}
\sum_{k=n+1}^{\infty} \frac{\cos(k\,t)}{k^{\nu}}(-1)^{k} = \sum_{k = n+1}^{\infty}\frac{\cos\left(k\,(t+\pi)\right)}{k^{\nu}}.
\end{equation}
is used. It follows from the above identity and the \cref{Coro:CoeffsPLUS}, the leading tern of pointwise error is
\begin{equation*}
f(x)-f_{n} \sim -\frac{\sin\left(\frac{2n+1}{2}(t+\pi)\right)}{2\sin(\frac{t+\pi}{2})}\,  \left(-\frac{2\,\beta_{1}}{\pi}\sin(\delta\pi)\,
\Gamma(2\delta+1)\,\ln(n) + \mu\,\beta_{1}\beta_{2}\right) \frac{\ln^{\mu-1}} {n^{2\delta +1 }}+
 \mathcal{O}\left( \frac{\ln^{\mu-2}(n)}{n^{2\delta+1}}\right).
\end{equation*}

Next, the item $(ii)$ is to be proved. For the case $x = -1$, the corresponding error of Chebyshev truncation is asymptotic to
\begin{equation*}
\begin{aligned}
f(x) - f_{n}&\sim  \sum_{k = n + 1}^{\infty} \left(-\frac{2\,\beta}{\pi}\right)\,\Gamma(2\delta + 1)\,\sin(\delta\,\pi)\, (-1)^{k}\,
\cos(k\,\pi)\,\frac{\ln^{\mu-1}(k)}{k^{2\delta +1}} \\
& \quad + \sum_{k = n + 1}^{\infty} \mu\,\beta_{1}\,\beta_{2} \frac{\ln^{\mu -1}(k)}{k^{2\delta + 1}}(-1)^{k}\cos(k\pi) +
\mathcal{O}\left(\frac{\ln^{m-2}(n)}{n^{2\delta}}\right).\\
\end{aligned}
\end{equation*}
By the \cref{Lem:PsiLogAsym}, the error can be asymptotically written as
\begin{equation*}
f(x) -f_{n} \sim -\frac{1}{\delta}\,\frac{\beta_{1}}{\pi}\,\Gamma(2\delta + 1) \, \sin(\delta\,\pi)\,\frac{\ln^{\mu}(n)}{n^{2\delta}} +
\frac{\mu}{2\delta}\,\beta_{1}\,\beta_{2}\, \ln^{\mu-1}{n^{2\delta}} + \mathcal{O}\left( \frac{ \ln^{\mu-2}(n) }{ n^{2\delta} } \right).
\end{equation*}
\end{proof}

\begin{figure}[!h]
\vspace{-0.3cm}
\includegraphics[scale=0.533]{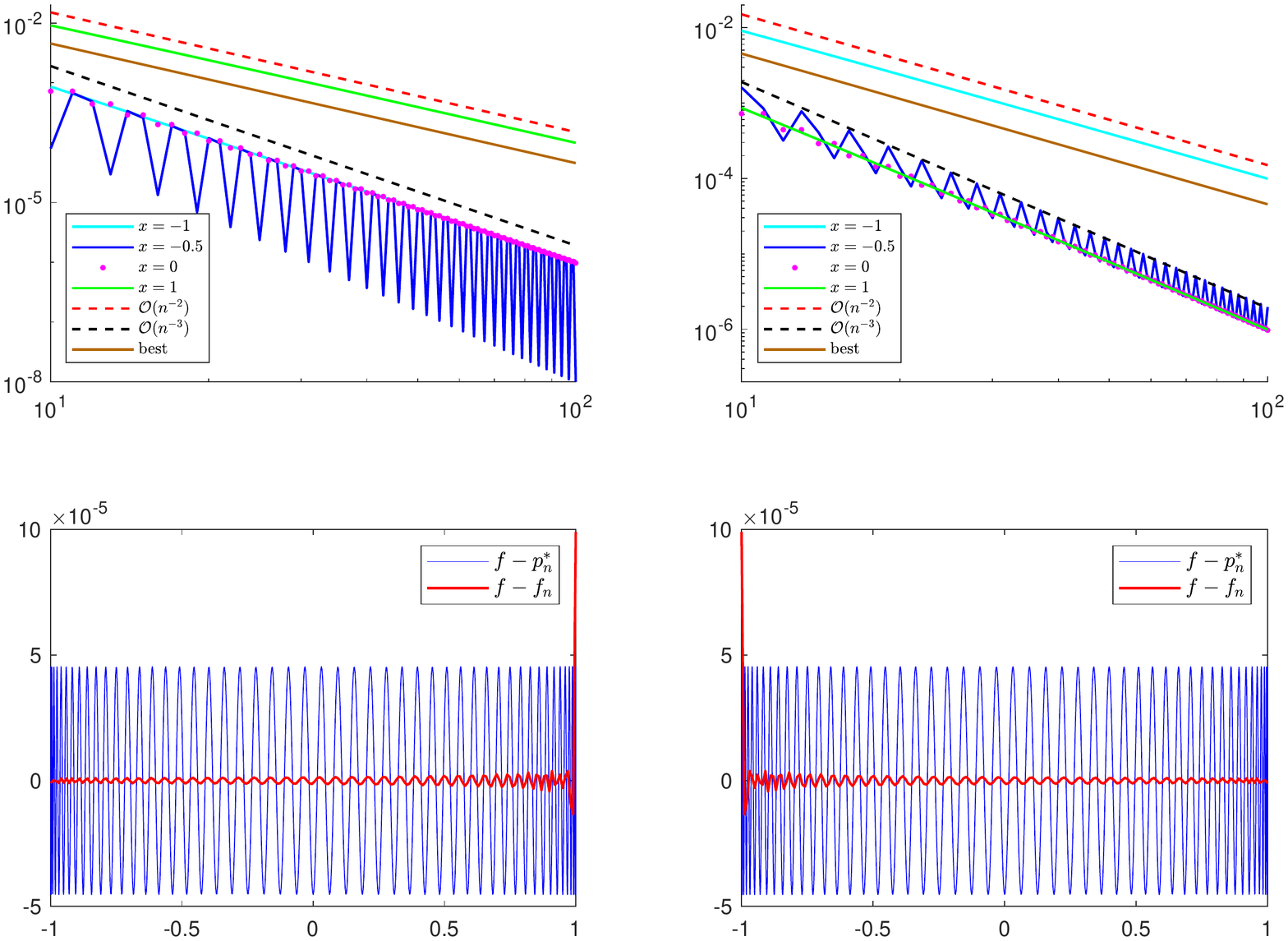}
\vspace{-1cm}
\caption{Pointwise error for logarithm singularity function. Upper left: Pointwise errors at $x=-1$ (Cyan line), $x=-0.5$(Blue line),
$x=0$(Magenta dot line), $x=1$(Green line) respectively for the function $f(x)=(1-x)\ln(1-x)$; Golden line : the error of the best approximation in the maximum norm for the same function; Red dash line : $\mathcal{O}(n^{-2})$; Black dash line : $\mathcal{O}(n^{-3})$. Upper right : Pointwise errors at the same points as ones on the left for the left boundary logarithm singularity function $f(x)=(1+x)\ln(1+x)$;
Golden line : the error of the best approximation in the maximum norm. Lower: the comparisons of the errors of Chebyshev truncations with
the ones of the best approximations. Lower left : the function is $f(x)=(1-x)\ln(1-x)$; Lower right : the function is $f(x)=(1+x)\ln(1+x)$. }
\label{fig:pointwiseErrlog}
\end{figure}

The theorem implies that the pointwise error at the singularity $x=1$ or $x= -1$ for the logarithm singularity function is $\mathcal{O}(n)$
larger than other points that is not singular, as is shown in \cref{fig:pointwiseErrlog}. Just as the error in $L_{\infty}$ norm on whole interval
is $O(n)$ larger than the interior interval that cut off the logarithm singularity boundary. We can also see that the
error of the best approximation and the error of the Chebyshev truncation decay at a same power of $n$, but the former is only a constant
multiplier less than the latter. Thus, the \cref{thm:ChebTrunNearBest} is so easy to make readers misunderstand that the Chebyshev truncation
is always inferior to the best approximation on the whole interval. Fortunately, Trefethen \cite{Trefethen2020} have noticed that it is not always right for the algebraic singularity function. Wang \cite{Wang2021best} have provided a firm theory to explain this phenomenon. Based on their great work, we make a further step for analysing the logarithm singularity function. In fact, the Chebyshev approximation is superior to the best approximation on almost whole interval except the narrow, narrow singular boundary layer for the logarithm singularity function. It is also easy to obtain that the leading term of error at the point $x= 1$ for the function \cref{eq:funoneminius} is same as the one at point $x=-1$ for the function \cref{eq:funtwoplus} if $g_{1}(1)= g_{2}(-1)$. Thus, if $g_{1}(x) = g_{2}(x)$, the errors of the both Chebyshev truncations in $L_{\infty}$ norm are same for corresponding functions.

Now we can explain why the pointwise error on the left boundary (at $x=-1$) is smaller than the pointwise error on the right boundary(at $x=1$)
for the function $f(x) = (1+\frac{x}{2})(1-x^2)\ln(1-x^2)$ as is shown in \cref{intro:AlgebraicVsLog}. At first glance, the result seems so surprising. Maybe some puzzles that the points $x=-1$ and $x=1$ are both singularity points. Why their behavior are different. In fact, the result is not surprising, the reason is following. For the function $f(x) = (1+\frac{x}{2})(1-x^2)\ln(1-x^2)$, we can write as
\begin{equation*}
\begin{aligned}
f(x) &= \left(1 + \frac{x}{2}\right)(1 + x)(1 - x)\ln(1 - x) + \left(1 + \frac{x}{2}\right)(1 - x)(1 + x)\ln(1 + x) \\
&= g_{1}(x)(1 - x)\ln(1 - x) + g_{2}(x)(1 + x)\ln(1 + x)
\end{aligned}
\end{equation*}
It is easy to obtain $g_{1}(1) = 3\, g_{2}(-1)$, thus the pointwise error at $x=1$ is three times as large as the error at $x=-1$.
It is caused by the function $g_{1}(x)$ and $g_{2}(x)$ rather than the singularities at $x=-1$ and $x=1$.
 Just as the description pointed out before, if $g_{1}(x)  = g_{2}(x)$, the maximum errors at two points $x=-1$ and $x=1$ are same.

Indeed, before this the authors \cite{zhang2021asymptotic} have noticed the maximum pointwise error ($L^{\infty}$ norm) of the truncated Chebyshev series is
$\mathcal{O}(n^{-2\gamma - 1})$ in the interior of the interval $[-1,1]$ , but the maximum pointwise error is $\mathcal{O}(n^{-2\gamma})$ in the
whole interval $[-1,1]$, for the function $f(x) = g(x)(1-x^2)^{\gamma}\ln(1-x^2)$.

\section{Pointwise error estimates of Chebyshev interpolation}\label{sec:interp}

In this section, the pointwise error of Chebyshev interpolations are analyzed. The two most widely used interpolants are considered.

\subsection{Chebyshev points of the first kind}
Let $\{x_{j}\}_{j=0}^{n}$ be the Chebyshev points of the first kind, i.e., $x_{j} = \cos(\frac{2j+1}{2n+2}\pi), j= 0,1,\cdots,n$ and
$p_{n}^{I}$ denotes the the polynomial of degree $n$ determined by the $n+1$ Chebyshev points of the first kind $\{x_{j}\}_{j=0}^{n}$. Based on
the discrete orthogonality, the Chebyshev interpolant can be represented as following
\begin{equation}
p_{n}^{I} = \sum_{k=0}^{n} {'} b_{k}\,T_{k}(x), \quad b_{k} = \frac{2}{n+1} \sum_{j=0}^{n} f(x_{j})\,T_{k}(x_{j}).
\end{equation}
The relationship of the truncation series coefficients and the interpolant coefficients for a finite $n$ is given
\begin{equation}
b_{k} = a_{k} + \sum_{\ell = 1}^{\infty}\,(-1)^{\ell}\, \left(a_{2\ell(n+1) - k} + a_{2\ell(n+1) + k}\right),\quad k = 0,1,\cdots, n.
\end{equation}
The detail derivation of the relation is given in \cite{FoxParker1968,Wang2021best}. The error of Chebyshev interpolant $p_{n}^{I}$ is
\begin{equation}\label{eq:errinpolantI}
f(x) - p_{n}^{I}(x) = \sum_{k = n+1}^{\infty} a_{k}\,T_{k}(x) + \sum_{k=0}^{n} \sum_{\ell = 1}^{\infty}\, (-1)^{\ell + 1}
\left(a_{2\ell(n+1) - k} + a_{2\ell(n+1) + k}\right)\,T_{k}(x).
\end{equation}
Next the second term of \cref{eq:errinpolantI} need to be estimated
\begin{equation}
\begin{aligned}
\sum_{k=0}^{n} \sum_{\ell = 1}^{\infty}(-1)^{\ell+1}\left(a_{2\ell(n+1) - k} + a_{2\ell(n+1) + k}\right)\,T_{k}(x) \sim \sum_{k=0}^{n}
\left(a_{2(n+1) -k} + a_{2(n+1) + k}\right) \cos(k\,t).\\
\end{aligned}
\end{equation}
Since
\begin{equation}
1 + \cos(t) + \cos(2\,t) + \cos(3\,t) + \cos(4\,t) + \cdots + \cos(n\,t) = \frac{\sin\left( (2n+1) \frac{t}{2} \right)+ \sin(\frac{t}{2})}{2\,
\sin(\frac{t}{2})}
\end{equation}
is uniformly bounded on $t\in [\delta_{1},\pi]$ for any small positive number $\delta_{1}$. Thus, one can further obtain
\begin{equation*}
 \left|\sum_{k=0}^{n} \left(a_{2(n+1) -k}
+ a_{2(n+1) + k}\right) \cos(k\,t)\right| \\
\le \frac{2}{\left|\sin(\frac{t}{2})\right|}\left(\left|a_{n+2}\right| + 3\left|a_{2n + 2}\right| + 2\left|a_{3n+2}\right|\right).
\end{equation*}
Finally, the interpolant error estimates are given in the following.

For the function \cref{eq:funoneminius}, one obtains the pointwise error is
\begin{equation}
\left|f(x) - p_{n}^{I}(x)\right| =
\begin{cases}
\mathcal{O}\left(\frac{\ln^{\mu}(n)}{n^{2\gamma + 1}}\right), & x\in[-1,1),\, \gamma\, \text{is not an integer},   \\
\mathcal{O}\left(\frac{\ln^{\mu-1}(n)}{n^{2\gamma + 1}}\right), & x\in[-1,1),\, \gamma\, \text{is an integer},\\
\mathcal{O}\left(\frac{\ln^{\mu}(n)}{n^{2\gamma}}\right), & x = 1,\,\gamma\, \text{is not an integer},\\
\mathcal{O}\left(\frac{\ln^{\mu-1}(n)}{n^{2\gamma}}\right), & x = 1,\,\gamma \,\text{is an integer}.
\end{cases}
\end{equation}
Similarly, for the function \cref{eq:funtwoplus}, one can also obtain that
\begin{equation}
\left|f(x) - p_{n}^{I}(x)\right| =
\begin{cases}
\mathcal{O}\left(\frac{\ln^{\mu}(n)}{n^{2\delta + 1}}\right), & x\in(-1,1],\, \delta\, \text{is not an integer},   \\
\mathcal{O}\left(\frac{\ln^{\mu-1}(n)}{n^{2\delta + 1}}\right), & x\in(-1,1],\, \delta\, \text{is an integer},\\
\mathcal{O}\left(\frac{\ln^{\mu}(n)}{n^{2\delta}}\right), & x = -1,\,\delta\, \text{is not an integer},\\
\mathcal{O}\left(\frac{\ln^{\mu-1}(n)}{n^{2\delta}}\right), & x = -1,\,\delta \,\text{is an integer}.
\end{cases}
\end{equation}

\subsection{Chebysehv points of the second kind}
Let $\{x_{j}\}_{j=0}^{n}$ be the set of the Chebyshev points of the second kind, i.e., $x_{j} = \cos\left(\frac{j}{n}\pi\right),\,j=0,1,\cdots,n$
and let $p_{n}^{II}$ be the interpolation polynomial of degree $n$ obtained by the Chebyshev points of the second. Based on the discrete
orthogonality of the Chebyshev polynomials, the Chebyshev interpolant $p_{n}^{II}$ can be written as
\begin{equation}
p_{n}^{II} = \sum_{k=0}^{n}{''}c_{k}T_{k}(x), \quad c_{k} = \sum_{j=0}^{n}{''}\frac{2}{n}f(x_{j})T_{k}(x_{j}).
\end{equation}
where the double prime denotes the first term and last term is to be halved. The relations between the coefficients $a_{k}$ of infinite series and
the coefficients $c_{k}$ of interpolant approximation used the Chebyshev points of the second kind is
\begin{equation}
c_{k} = a_{k} + \sum_{\ell = 1}^{\infty} \left(a_{2\ell n-k}+ a_{2\ell n+k}\right).
\end{equation}
The error is
\begin{equation}\label{eq:errChebSecond}
f(x) - p_{n}^{II} = \sum_{k = n+1}^{\infty}a_{k}T_{k}(x) - \sum_{k=0}^{n} \sum_{\ell = 1}^{\infty} \left(a_{2\ell n-k}+ a_{2\ell n+k}\right)
T_{k}(x).
\end{equation}
To analyze the error of the interpolant approximation based on the Chebysehv points of second kind, the second term \cref{eq:errChebSecond} need to
be estimated. Similar analysis to one for the Chebyshev points of first kind, the errors of interpolant for the functions \cref{eq:funoneminius}
, \cref{eq:funtwoplus} are given in the following.

For the function \cref{eq:funoneminius}, one obtains the pointwise error is
\begin{equation}
\left|f(x) - p_{n}^{II}(x)\right| =
\begin{cases}
\mathcal{O}\left(\frac{\ln^{\mu}(n)}{n^{2\gamma + 1}}\right), & x\in[-1,1),\, \gamma\, \text{is not an integer},   \\
\mathcal{O}\left(\frac{\ln^{\mu-1}(n)}{n^{2\gamma + 1}}\right), & x\in[-1,1),\, \gamma\, \text{is an integer},\\
\mathcal{O}\left(\frac{\ln^{\mu}(n)}{n^{2\gamma}}\right), & x = 1,\,\gamma\, \text{is not an integer},\\
\mathcal{O}\left(\frac{\ln^{\mu-1}(n)}{n^{2\gamma}}\right), & x = 1,\,\gamma \,\text{is an integer}.
\end{cases}
\end{equation}

For the function \cref{eq:funtwoplus}, one can also obtain that
\begin{equation}
\left|f(x) - p_{n}^{II}(x)\right| =
\begin{cases}
\mathcal{O}\left(\frac{\ln^{\mu}(n)}{n^{2\delta + 1}}\right), & x\in(-1,1],\, \delta\, \text{is not an integer},   \\
\mathcal{O}\left(\frac{\ln^{\mu-1}(n)}{n^{2\delta + 1}}\right), & x\in(-1,1],\, \delta\, \text{is an integer},\\
\mathcal{O}\left(\frac{\ln^{\mu}(n)}{n^{2\delta}}\right), & x = -1,\,\delta\, \text{is not an integer},\\
\mathcal{O}\left(\frac{\ln^{\mu-1}(n)}{n^{2\delta}}\right), & x = -1,\,\delta \,\text{is an integer}.
\end{cases}
\end{equation}

\begin{figure}[!ht]
\vspace{-0.5cm}
\includegraphics[height=2in,width=6.2in]{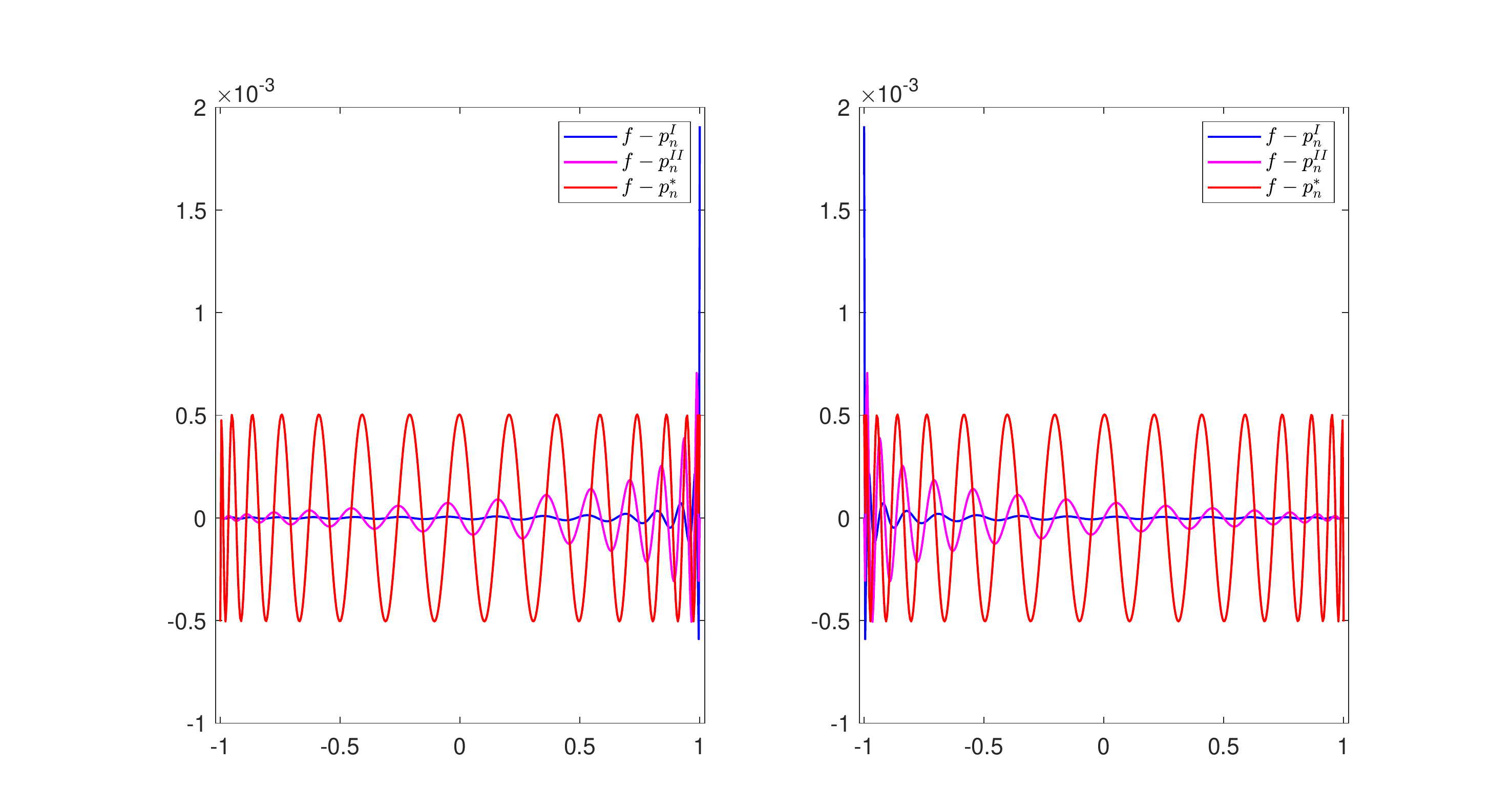}
\caption{Pointwise errors of $p_{n}^{I}$, $p_{n}^{II}$ and $p_{n}^{*}$ for the function $f(x) = (1-x)\ln(1-x)$(Left) and $f(x) = (1+x)\ln(1+x)$(Right) with
$n= 30$.}\label{fig:PointwiseErrorInterplants}
\end{figure}

It deserves to point out the pointwise errors are not equal of the two kind interpolant approximations for a logarithm singularity function.
Although we can not give a define conclusion which is better, the errors at the weak singularity points are considered in the following. Here the
function \cref{eq:funoneminius} is considered, other logarithm singularity functions can be analysed in the similar way. For the
function \cref{eq:funoneminius}, the weak singular point is $x=1$. The error of the interpolant $p_{n}^{I}$ at this point is
\begin{equation}\label{eq:errFirstInpAtone}
\left.[f - p_n^{I}]\right|_{x=1} = \sum_{k=n+1}^{\infty}a_{k} - \sum_{k=0}^{n} \sum_{\ell = 1}^{\infty}\, (-1)^{\ell + 1}
\left(a_{2\ell(n+1) - k} + a_{2\ell(n+1) + k}\right),
\end{equation}
and the error the interpolant $p_{n}^{II}$ at $x=1$ is
\begin{equation}\label{eq:errSecondInpAtone}
\left.[f - p_n^{II}]\right|_{x=1} = \sum_{k=n+1}^{\infty}a_{k} - \sum_{k=0}^{n} \sum_{\ell = 1}^{\infty}\,
\left(a_{2\ell n - k} + a_{2\ell n + k}\right).
\end{equation}
By the \cref{thm:asycoeffsfunminus}, the sign of the coefficients $a_{n}$ will not change when $n$ is more than a definite positive integer. To simplify
the analysis, without loss of the generality we assume the coefficients $a_{n}$ are positive when $n$ is more than a definite positive integer. The first
term on the right in \cref{eq:errFirstInpAtone} is same as one in \cref{eq:errSecondInpAtone}. It is not difficult to see that the second term
in \cref{eq:errSecondInpAtone} is much larger than the one in \cref{eq:errFirstInpAtone}. Thus the absolute value of error of interpolant $p_n^{I}$ is larger than one of interpolant $p_{n}^{II}$, as is shown in \cref{fig:PointwiseErrorInterplants}.

\section{Conclusions}
In \cite{Xiang2021}, the author gave the asymptotic Chebyshev coefficients for the functions with logarithms regularities based on the Hib type formula between Jocobi polynomials and Bessel functions. Different with the proposed methods, the steepest descent method is used to evaluate the asymptotic Chebyshev coefficients for the infinite series of a function with logarithm regularities \cref{eq:orignalfun}. Then the pointwise error of Chebyshev truncations are estimated. For the function \cref{eq:orignalfun}, the pointwise errors at the singularity boundaries are $\mathcal{O}(n)$ larger than ones at the interior points. This provides a more through explaination the error behavior of the Chebyshev truncations, furthermore, it also helps us make a full understand the real relation between the best polynomial approximation and the Chebshev truncation. In the last section showed that the behaviors of pointwise errors of interpolant approximations based on the two kind interpolation points are analogous to the ones of the Chebyshev truncation errors.

Finally, we echo a similar sentence in \cite{Trefethen2020}, best polynomial approximation is optimal in uniform norm, but Chebyshev interpolants are better for the functions with logarithmic regularities.

\section*{Acknowledgments}

This work supported by the National Natural Science Foundation of China (No.12101229) and the Hunan Provincial Natural Science Foundation of China (No.2021JJ40331).

\bibliographystyle{siamplain}
\bibliography{references}
\end{document}

%% file: ex_shared.tex
\usepackage{lipsum}
\usepackage{amsfonts}
\usepackage{graphicx}
\usepackage{epstopdf}
\usepackage{algorithmic}
\ifpdf
  \DeclareGraphicsExtensions{.eps,.pdf,.png,.jpg}
\else
  \DeclareGraphicsExtensions{.eps}
\fi

\newsiamremark{remark}{Remark}
\newsiamremark{hypothesis}{Hypothesis}
\crefname{hypothesis}{Hypothesis}{Hypotheses}
\newsiamthm{claim}{Claim}

\headers{Mysteries of the Chebyshev expansion}{X. Zhang}

\title{The mysteries of the best approximation and Chebyshev expansion for the function with logarithmatic regularities}
\author{Xiaolong Zhang\thanks{Key Laboratory of Computing and Stochastic Mathematics (Ministry of Education), School of Mathematics and Statistics, Hunan Normal University, Changsha, 410081, China
  (\email{xlzhang@hunnu.edu.cn}).}
}

\usepackage{amsopn}
